\documentclass[11pt]{amsart}
\usepackage{amscd,amsmath,amssymb,amsfonts}
\usepackage{color}
\usepackage[cmtip, all]{xy}
\usepackage{graphicx}
\usepackage{enumerate}
\usepackage{appendix}
\numberwithin{equation}{section}
\usepackage[usenames,dvipsnames]{xcolor}
\usepackage{times}
\usepackage{amsxtra}
\usepackage{amsthm}
\usepackage{footnote}
\usepackage{footmisc}

\newcounter{savefootnote}
\newcounter{symfootnote}
\newcommand{\symfootnote}[1]{%
   \setcounter{savefootnote}{\value{footnote}}%
   \setcounter{footnote}{\value{symfootnote}}%
   \ifnum\value{footnote}>8\setcounter{footnote}{0}\fi%
   \let\oldthefootnote=\thefootnote%
   \renewcommand{\thefootnote}{\fnsymbol{footnote}}%
   \footnote{#1}%
   \let\thefootnote=\oldthefootnote%
   \setcounter{symfootnote}{\value{footnote}}%
   \setcounter{footnote}{\value{savefootnote}}%
}

\newcommand{\nc}{\newcommand}

\numberwithin{equation}{section}

\swapnumbers

\newtheorem{problem*}{Problem}[section]
\newtheorem{key idea}[thm]{Key idea}
\newtheorem*{theorem*}{Theorem}
\newtheorem*{theoremT*}{Taniyama's ``theorem"}
\newtheorem*{conjI*}{Conjecture I}
\newtheorem*{conjII*}{Conjecture II}
\newtheorem*{conj*}{Conjecture}
  \newtheorem{innercustomgeneric}{\customgenericname}
\providecommand{\customgenericname}{}
\newcommand{\newcustomtheorem}[2]{%
  \newenvironment{#1}[1]
  {%
   \renewcommand\customgenericname{#2}%
   \renewcommand\theinnercustomgeneric{##1}%
   \innercustomgeneric
  }
  {\endinnercustomgeneric}
}

\newcustomtheorem{customthm}{Theorem}
\newcustomtheorem{customconj}{Conjecture}
\newcustomtheorem{customconjI}{Conjecture I}
\newcustomtheorem{customconjII}{Conjecture II}

\nc{\cf}{{\mathcal F}}
\nc{\la}{\lambda}
\nc{\ct}{T^*X}
\nc{\str}{{\mathcal O}_X(\lambda)}
\nc{\od}{\operatorname  {D}}
\nc{\occ}{\operatorname {CC}}
\nc{\ode}{\operatorname {E}}
\nc{\of}{\operatorname {F}}
\nc{\op}{\operatorname {P}}
\nc{\om}{\operatorname {M}}
\nc{\pro}{G_{\bold R}\times H \times \hat X}
\nc{\Ext}{\operatorname{Ext}}
\nc{\Hom}{\operatorname{H}} 
\nc{\Char}{\operatorname{Char}}

\nc{\Bun} {{\operatorname{Bun}}}
\nc {\Aut }{{\operatorname{Aut}}}
\nc{\V}{{\operatorname {Vec}_\Bbbk}}
\nc{\oh}{{\operatorname H}}
\nc{\cod}{{\operatornamewithlimits{codim}}}
\nc{\dbs}{\operatorname {D}_{\cs}^b}
\nc{\ps}{\operatorname {P}_{\cs}}
\nc{\p} {\operatorname {P}}

\nc{\GO}{{G(\mathcal O)}}
\nc{\GK}{{G(\mathcal K)}}
\nc{\G}{{\mathcal G}}
\nc{\fa}{{\mathfrak a}}
\nc{\fb}{{\mathfrak b}}
\nc{\fc}{{\mathfrak c}}
\nc{\fd}{{\mathfrak d}}
\nc{\fe}{{\mathfrak e}}
\nc{\ff}{{\mathfrak f}}
\nc{\fg}{{\mathfrak g}}
\nc{\fh}{{\mathfrak h}}
\nc{\fiI}{{\mathfrak i}}  
\nc{\ffi}{{\mathfrak i}}  
\nc{\fj}{{\mathfrak j}}
\nc{\fk}{{\mathfrak k}}
\nc{\fl}{{\mathfrak{l}}}
\nc{\fm}{{\mathfrak m}}
\nc{\fn}{{\mathfrak n}}
\nc{\fo}{{\mathfrak o}}
\nc{\fp}{{\mathfrak p}}
\nc{\fq}{{\mathfrak q}}
\nc{\fr}{{\mathfrak r}}
\nc{\fs}{{\mathfrak s}}
\nc{\ft}{{\mathfrak t}}
\nc{\fu}{{\mathfrak u}}
\nc{\fv}{{\mathfrak v}}
\nc{\fw}{{\mathfrak w}}
\nc{\fz}{{\mathfrak z}}
\nc{\fx}{{\mathfrak x}}
\nc{\fy}{{\mathfrak y}}

\nc{\fS}{{\mathfrak S}}

\nc{\cA}{{\mathcal A}}
\nc{\cB}{{\mathcal B}}
\nc{\cC}{{\mathcal C}}
\nc{\cD}{{\mathcal D}}
\nc{\cE}{{\mathcal E}}
\nc{\cF}{{\mathcal F}}
\nc{\cG}{{\mathcal G}}
\nc{\cH}{{\mathcal H}}
\nc{\cI}{{\mathcal I}}
\nc{\cJ}{{\mathcal J}}
\nc{\cK}{{\mathcal K}}
\nc{\cL}{{\mathcal L}}
\nc{\cM}{{\mathcal M}}
\nc{\cN}{{\mathcal N}}
\nc{\cO}{{\mathcal O}}
\nc{\cP}{{\mathcal P}}
\nc{\cQ}{{\mathcal Q}}
\nc{\cR}{{\mathcal R}}
\nc{\cS}{{\mathcal S}}
\nc{\cT}{{\mathcal T}}
\nc{\cU}{{\mathcal U}}
\nc{\cV}{{\mathcal V}}
\nc{\cW}{{\mathcal W}}
\nc{\cZ}{{\mathcal Z}}
\nc{\cX}{{\mathcal X}}
\nc{\cY}{{\mathcal Y}}

\nc{\bA}{{\mathbb A}}
\nc{\bB}{{\mathbb B}}
\nc{\bC}{{\mathbb C}}
\nc{\bD}{{\mathbb D}}
\nc{\bE}{{\mathbb E}}
\nc{\bF}{{\mathbb F}}
\nc{\bG}{{\mathbb G}}
\nc{\bH}{{\mathbb H}}
\nc{\bI}{{\mathbb I}}
\nc{\bJ}{{\mathbb J}}
\nc{\bK}{{\mathbb K}}
\nc{\bL}{{\mathbb L}}
\nc{\bM}{{\mathbb M}}
\nc{\bN}{{\mathbb N}}
\nc{\bO}{{\mathbb O}}
\nc{\bP}{{\mathbb P}}
\nc{\bQ}{{\mathbb Q}}
\nc{\bR}{{\mathbb R}}
\nc{\bS}{{\mathbb S}}
\nc{\bT}{{\mathbb T}}
\nc{\bU}{{\mathbb U}}
\nc{\bV}{{\mathbb V}}
\nc{\bW}{{\mathbb W}}
\nc{\bZ}{{\mathbb Z}}
\nc{\bX}{{\mathbb X}}
\nc{\bY}{{\mathbb Y}}

\nc{\La}{{\Lambda}}

\newcommand{\QQ}{{\mathbb Q}}

\newcommand{\ZZ}{{\mathbb Z}}

\nc{\bi}{\bibitem}

\renewenvironment{thebibliography}[1]{%
\begin{oldthebibliography}{#1}%
\setlength{\parskip}{0ex}%
\setlength{\itemsep}{0ex}%
}%
{%
\end{oldthebibliography}%
}

\nc{\Ker}{\operatorname{Ker}}
\nc{\Coker}{\\operatorname{Coker}}

\nc{\codim}{\operatorname{codim}}

\begin{document} 
\author{Michael Harris}
\title[VIRTUES OF PRIORITY]{Virtues of Priority}

\address{Michael Harris\\
Department of Mathematics, Columbia University, New York, NY  10027, USA}
 \email{harris@math.columbia.edu}

\dedicatory{In memory of Serge Lang}
\thanks{This article was originally written in response to an invitation by a group of philosophers as part of ``a proposal for a special issue of the philosophy journal {\it Synth\`ese} on virtues and mathematics." The invitation read, ``We would be delighted to be able to list you as a prospective contributor. This would, of course, be in no way binding on either party." The author accepted the invitation on the assumption that the account it contains of the development of a particularly bitter controversy within the mathematical community might prove useful as raw material for philosophers of mathematics. Unfortunately, the referees {\it Synth\`ese} chose recommended that the article be entirely rewritten to respect the professional norms of philosophy of mathematics. Not being a philosopher of mathematics, the author had no intention of doing this, and the editors of {\it Synth\`ese} saved the author the trouble by rejecting the article.
~~~~~~~~~~~~~~~~~~~~~~~~~~~~~This work was partially supported by NSF Grant DMS-1701651.  This work was also supported by the National Science Foundation under Grant No. DMS-1440140 while the author was in residence at the Mathematical Sciences Research Institute in Berkeley, California, during the Spring 2019 semester.}

\maketitle

\section*{Introduction:  Originality and other virtues}

If hiring committees are arbiters of mathematical virtue, then letters of recommendation should give a good sense of the virtues most appreciated by mathematicians.  You will not see ``proves true theorems" among them.  That's merely part of the job description, and drawing attention to it would be analogous to saying an electrician won't burn your house down, or a banker won't steal from your account.  I don't know how electricians or bankers recommend themselves to one another, but I have read a lot of letters for jobs and prizes in mathematics, and their language is revealing in its repetitiveness.  Words like ``innovative" or ``original" are good, ``influential" or ``transformative" are better, and ``breakthrough" or ``decisive" carry more weight than ``one of the best."  Best, of course, is ``the best," but it is only convincing when accompanied by some evidence of innovation or influence or decisiveness.  

When we try to answer the questions: what is being innovated or decided? who is being influenced?  -- we conclude that the virtues highlighted in reference letters point to mathematics as an undertaking relative to and within a community.  This is hardly surprising, because those who write and read these letters do so in their capacity as representative members of this very community.  Or I should say: members of overlapping communities, because the virtues of a branch of mathematics whose aims are defined by precisely formulated conjectures (like much of my own field of algebraic number theory) are very different from the virtues of an area that grows largely by exploring new phenomena in the hope of discovering simple underlying principles (like much of the field of dynamical systems).  The interpretation of letters is left to the experts who apply not only their knowledge of their field's problems but also their familiarity with its value system.  The same vocabulary (original, influential, breakthrough) can thus apply to developments whose character differs radically from one field to another.  Some of the five mathematicians who in 2014 were the first laureats of Silicon Valley's {\it Breakthrough Prizes} were honored for solving longstanding unsolved problems, while others defined new classes of objects or raised original questions.  

``Originality" is a virtue recognized across disciplinary boundaries.  The ability to understand and explain the details of another mathematician's difficult proof is undoubtedly admired as a virtue, but it's almost never enough to get you a job.  For that you have to have brought something {\it new} to your field.  This sounds so banal to be unworthy of attention, but the philosophical problem is glaringly obvious:  how do we determine what is and is not {\it new}?  Here the word {\it new} is used as shorthand for the more complex notion of ``bringing something new to your field."  Bare newness is not enough.  You can take 99 theorems in your field and string them together into a single sentence connected by ``and"s but that's not {\it new} even if no one has ever thought to do it before.  As I already hinted, most true theorems are not {\it new} in the virtuous sense I share with my fellow members of the hiring or prize committee.  

Surely an element of taste is involved.  It would be easiest to leave the clarification of the virtue of newness to specialists in mathematical aesthetics, who are in as difficult a position as any philosophers of mathematics.  Specifically, they face the difficulty common to all philosophers who deal with collective characteristics of communities:  how to distinguish the abstract (and presumably permanent) principle they wish to emphasize from the (presumably transitory) idiosyncracies of individual members of the community?  Or more crudely, is a mathematical contribution {\it innovative} because it is claimed to be so by influential or powerful members of the community -- the kind who are asked to write letters of recommendation?  Are reference letters just politics by other means?   

These considerations can be left in the background when it has somehow come to be accepted that answering a certain question -- solving a famous conjecture, for example -- meets the community's standards of importance.  No member in good standing of the community of geometers doubted that Grigori Perelman's claim to have solved the Poincar\'e conjecture was {\it new} in the relevant sense, because (a) the problem had long been considered central to the field; (b) no one had solved it before Perelman; and (c) several teams of experts worked through the proof and vouched for its correctness.  

It is already more difficult to characterize the newness of a proof of a theorem that was not long-awaited.  Judgments in such cases are obtained as the result of a process of implicit negotiation among those recognized as leaders within the community. On the basis of the standards invoked in such negotiations it should be possible to form a reasonably coherent (or perhaps chaotic) picture of the community's virtues.

Still more problematic are the mathematical ideas that express themselves in the form of {\it conjectures}.  It is at least conventionally believed that the {\it goal} of mathematics is to produce proofs of theorems.\footnote{However, William Thurston famously wrote that ``The product of mathematics is clarity and understanding. Not theorems, by themselves."  https://mathoverflow.net/posts/44213/revisions.}  So faced with such a proof, one at least knows where  the negotiation starts.  The place of conjectures is much more ambiguous:  they can guide the thoughts of mathematicians for generations, or even for centuries, but in principle a conjecture that turns out to be mistaken can also vanish without a trace.\footnote{The canonical reference for the role of conjecture within mathematical practice is still Barry Mazur's article on the subject:  Mazur, B., Conjecture. {\it Synth\`ese}, vol. 111, no. 2, 1997, pp. 197--210. JSTOR, www.jstor.org/stable/20117628.}    

\section*{The Modularity Conjecture:  Weil's 1967 paper}

In 1967 Andr\'e Weil published an 8-page article entitled ``\"Uber die Bestimmung Dirichletscher Reihen durch Funktionalgleichungen" -- ``On the determination of Dirichlet series by means of functional equations."  The title was an intentional allusion to a 1936 paper of Erich Hecke --``\"Uber  die Bestimmung Dirichletscher Reihen durch ihre Funktionalgleichung":  only one functional equation is used in Hecke's title and in his paper for determination, whereas the result of Weil's paper requires a collection of functional equations.  The main theorem of the paper, now known as {\it Weil's converse theorem}, provides conditions that guarantee that a Dirichlet series -- an infinite series of the form  $\sum_n a_n n^{-s}$, where $s$ is a complex variable -- arises by a well-known procedure from a {\it modular form} of a specified type.  I will return to this last point later.  

Weil must have considered his result respectable, because it was his contribution to a special issue of the journal {\it Mathematische Annalen} in honor of the 70th birthday of Carl Ludwig Siegel, a mathematician whom Weil greatly admired.  The last page or so of his brief text was devoted to ``reflections on the zeta functions of elliptic curves that perhaps also merit some attention" (berlegungen \"uber Zetafunktionen elliptischer Kurven ankn\"upfen, die vielleicht auch einige Aufmerksamkeit verdienen).  These reflections came very close to, but stopped short of proposing the conjecture, now solved, whose name has since been a matter of controversy:  

\begin{conjI*}  {\bf The zeta function of any elliptic curve with rational coefficients coincides with the zeta function of a modular form of the type specified in the main theorem of Weil's paper.  }
\end{conjI*}

I choose to formulate this sentence as a conjecture for the sake of clarity, and to set the conjecture in boldface type, but Weil presented his ``reflections" in normal type. Moreover, instead of coming out and stating a conjecture, he explained that this is sometimes the case -- here he cited a result of Goro Shimura -- but added that ``whether things always behave that way" appears for the moment to be problematic (scheint im Moment noch problematisch zu sein) and may be recommended to the interested reader as an exercise (mag dem interessierten Leser als \"Ubungsaufgabe empfohlen werden).  

These last sentences came at the end of a paragraph that began with Weil's very strong hints in favor of believing Conjecture I:  ``On certain theoretical grounds one may with reasonable certainty suspect"\footnote{In the original German:  ``Aus gewissen theoretischen Gr\"unden darf man mit ziemlicher Sicherheit vermuten."  The German admits more than one translation.  The online dictionary https://www.linguee.com/german-english/translation/ provides sentences in which {\it mit ziemlicher Sicherheit} means ``almost certainly." Moreover, {\it Vermutung} is also the first and best-known German equivalent of the English word {\it conjecture,} which leads one to wonder whether Weil was being entirely sincere when in 1979 he argued against conjectures in general.}  that something is true that entails the claim, in view of the main result of the paper.    This ``suspicion" -- in a rather stronger form, mentioned in that same paragraph of Weil's paper -- soon was circulating under the name of the {\it Weil Conjecture}.  I will have to provide some background in order to explain the nature of these hints, but for the moment let me focus on the superficial aspect of the situation.  Weil has drawn his readers' attention to a series of steps, some of them proved in his paper, some ``reasonably certain,"  leading to the suspicion and, with another step that is ``of course obvious to expect" (nat\"urlich naheliegend zu erwarten), to the stronger form that will be explained in a moment.

\section*{Lightning introduction to elliptic curves, modular forms, and their Dirichlet series}

But first I had better review the mathematical notions invoked in the ``suspicion."  An {\it elliptic curve with rational coefficients} is a geometric object $E$ with number theoretic meaning:  it is the set of solutions to a polynomial equation in two variables of the form

\begin{equation}\label{ell}  y^2 = x^3 + ax + b  \tag{*}
\end{equation}			 

where $a$ and $b$ are rational numbers (with the property that the polynomial on the right-hand side of the equation has no multiple roots).    A modular form is a different kind  of geometric object with a very different number-theoretic meaning.  To explain this, we consider the set of complex numbers 
$\mathfrak{H} = \{z = x + iy, y > 0\}$-- the upper half-plane in the complex plane.  The $2 \times 2$ matrices 
$$ \begin{pmatrix} a & b \\ c & d \end{pmatrix}$$
  with determinant $1$ form a group $SL(2,\ZZ)$ that acts on $\mathfrak{H}$ by the formula
\begin{equation*} z \mapsto \frac{az+b}{cz+d}
\end{equation*}
 
Whenever $\Gamma$ is a subgroup of finite index in $SL(2,\ZZ)$, the quotient $\Gamma\backslash \mathfrak{H}$ can be identified with a compact Riemann surface with some points missing; adding the points yields the {\it modular curve} $X(\Gamma)$.  A modular form of weight $2$ and level $\Gamma$ is a complex analytic function $f(z)$ on $\mathfrak{H}$ that is not quite invariant under $\Gamma$; rather the differential form $f(z)dz$ is invariant under $\Gamma$.  

There is a procedure A for associating a Dirichlet series to such an $f$, and there is also a procedure B for associating a Dirichlet series -- Weil denotes it $\Lambda$ -- to an elliptic curve $E$.  Weil proves that any Dirichlet series that satisfies a property (A3) -- the ``functional equations" (plural) of Weil's title -- arises from procedure A -- with specified {\it level subgroup} $\Gamma$ -- and alludes to ``certain theoretical grounds" for believing that the Dirichlet series arising from procedure B also satisfy property (A3).  This motivates the original suspicion, which I now present in full:
\medskip

{\it Aus gewissen theoretischen Gr\"unden darf man mit ziemlicher Sicherheit vermuten, dass $\Lambda$ eine Funktionalgleichung besitzt... Weiter l\"asst sich vermuten} [that the functional equations of condition (A3) also hold upon twist by characters].  
\medskip

On the other hand, $E$ can also be viewed as a Riemann surface, and there is a stronger form of the conjecture:

\begin{conjII*}   {\bf Any elliptic curve $E$ with rational coefficients admits a modular parametrization:   a non-trivial map of Riemann surfaces from $X(\Gamma)$ to $E$, again with specified $\Gamma$.  }
\end{conjII*}

For nearly 20 years, such a map was called a {\it Weil uniformization}, in conferences and seminars (including many that I attended) as well as in a number of influential papers, reflecting the consensus that Weil had formulated an idea whose implications included the existence of such a map.  But then things got more complicated.

\section*{One conjecture, many names}

Actually, things were already more complicated within a few years of the publication of Weil's papers, because it turned out that Yutaka Taniyama had listed something very like Conjectures I and II as items 12 and 13 of a collection of 36 problems.  This list was published in Japanese in his collected works, but it was handed out in English at the 1955 Tokyo-Nikko conference on number theory, attended by Weil, Shimura, and Jean-Pierre Serre,\footnote{I am following Serge Lang's account in ``Some history of the Shimura-Taniyama Conjecture," {\it Notices of the AMS}, {\bf 42}, November 1995, 1301--1307. }  who also has a part to play in this story.  Great things were expected of Taniyama, but he committed suicide a few years after 1955 and is now mainly remembered for the theory of complex multiplication that he developed with Shimura, and for his part in formulating the boldface statement Conjecture I, whose official name had become the {\it Taniyama-Weil Conjecture}, no later than 1977.   By then it was generally understood that the statement would follow from the {\it Langlands program}, a vast and ambitious network of conjectures and theorems that linked {\it automorphic forms} -- generalizations of modular forms -- to generalizations of the number-theoretic structures derived from elliptic curves.  And it was generally believed that the Langlands program was a distant fantasy, in large part because no one had any idea how to begin to prove the (then-so-called) Taniyama-Weil Conjecture -- which I will henceforth call the {\it Modularity Conjecture}.

Positions hardened about halfway between the announcement by Frey and Ribet that Fermat's Last Theorem would be a consequence of the Modularity Conjecture, and the announcement by Wiles of his proof (later completed in a paper of Taylor-Wiles) of enough cases of the Modularity Conjecture to derive Fermat's Last Theorem.  In February 1988, before Wiles announced his result, Joseph Oesterl\'e had published a Bourbaki seminar talk whose first theorem, attributed to Ken Ribet, asserts that ``La conjecture de Taniyama-Weil implique le th\'eor\`eme de Fermat."   No paper by Shimura was cited in Oesterl\'e's bibliography.  Ribet himself, writing in the United States two years later, published a paper in France with the title ``From the Taniyama-Shimura conjecture to Fermat's last theorem," in which his reference to Weil's 1967 paper draws attention to what he considered to be Weil's contribution to the story.
\medskip

{\it Finally, A. Weil proved in [35] that an elliptic curve $E$ over $\QQ$ is modular provided that its $L$-function $L(E,s)$ has the analytic properties one expects of it.}
\medskip

Mathematicians' choice of terminology for the conjecture largely broke down along national lines, and remains contested nearly 20 years after a complete proof of the Modularity Conjecture was published by Christophe Breuil, Brian Conrad, Fred Diamond, and Richard Taylor.\footnote{Wikipedia is an international battleground.  The English, German, French, Russian, and Arabic pages are entitled ``The Modularity Theorem" but disagree on its former terminology as conjecture:  in English and German it was the ``Taniyama-Shimura Conjecture," in Russian the ``Taniyama-Shimura-Weil," and in French either ``Taniyama-Weil," ``Shimura-Taniyama-Weil," or ``Shimura-Taniyama."  (The Arabic page mentions Taniyama-Shimura-Weil ``and other names.") The Spanish, Italian, Swedish, Hebrew, Catalan (and, as far as I can tell, the Chinese) pages are entitled ``The Taniyama-Shimura Theorem"; the Portuguese page call it the ``Shimura-Taniyama-Weil Theorem"; and the Dutch page is ``The Shimura-Taniyama Theorem."  Wikipedia pages on mathematical topics are generally maintained by mathematicians, who, as far as I know, do not report to their national authorities.
The proof by Breuil et al. is in the article ``On the modularity of elliptic curves over $\QQ$: wild $3$-adic exercises", {\it Journal of the American Mathematical Society}, 14 (4): 843--939 (2001).}    Public disputes over priority were acrimonious on occasion, private disputes much more often.  My own discomfort with the situation, and my dissatisfaction with the quality of the arguments put forward in support of one position or the other, led me to seek clarity and consolation in philosophy.  Thus, by a circuitous route, I find myself writing an article for a philosophical journal\symfootnote{The article was written in response to an invitation but was ultimately rejected  by the journal, for reasons to be discussed elsewhere.} on the topic that revived my interest in philosophy in the first place.

\section*{One conjecture, many virtues:  comparing the cases for Taniyama, Shimura, and Weil}

In fact, although the difference of opinion had all the characteristics of a vulgar personality conflict -- it was common knowledge that relations between Serre and Shimura were difficult, and the pages Shimura devotes to Serre in his memoirs are so gratuitously nasty that one wonders why Springer-Verlag agreed to publish them -- arguments were advanced in support of each of the attributions.  It is these arguments that deserve to be evaluated for their implicit or explicit assumptions about the virtues attributed to priority.

\subsection*{1. The case for Taniyama}  By all accounts, Taniyama was the first to have expressed an expectation along the lines of the conjecture, in items 12 and 13 of his list of problems at the Tokyo-Nikko conference.  It is also agreed that his version was not literally correct as stated.\footnote{This opinion is shared by Lang, Serre, and Shimura.  In a letter to Lang, Shimura wrote that Taniyama ``was not completely careful, and if someone had pointed out this, he would have agreed that the problem would have to be revised accordingly."}    I copy the English translation of Taniyama's questions from Lang's article in the {\it AMS Notices}:

\hspace*{2.8 cm}\includegraphics{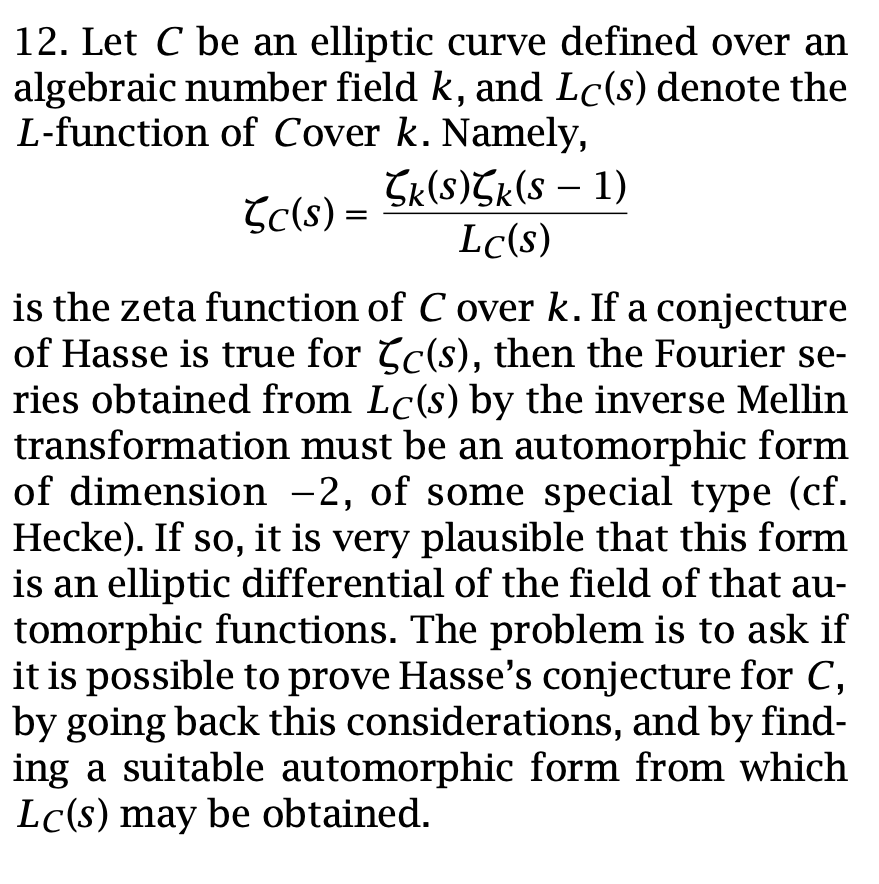}

\hspace*{2.8 cm}\includegraphics{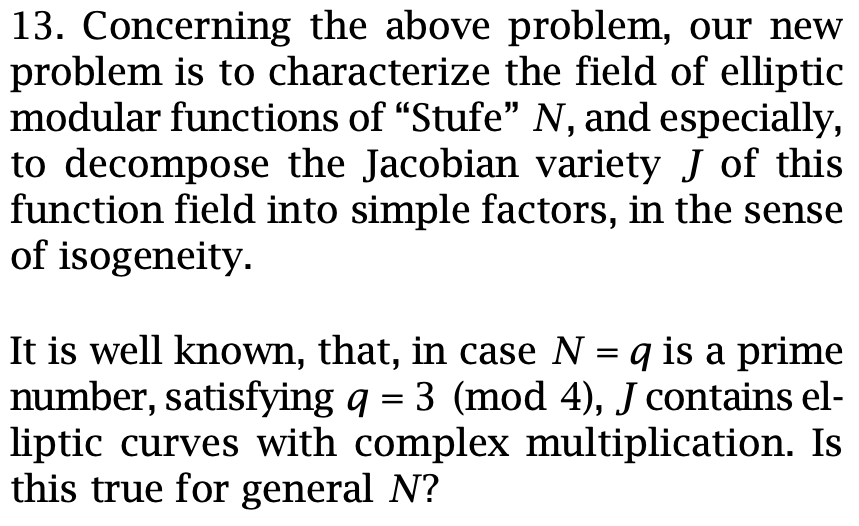}

\subsection*{2. The case for Shimura}  (a) The case for Shimura's priority was made by Lang in his 1995 article:  Lang established that Shimura had actually stated a version of the conjecture in the presence of Serre and Weil, at a specific place (the Institute for Advanced Study in Princeton) and during a specific time interval (``most likely in 1964" , according to a private message\footnote{Shimura, private emails to the author, May 6-12, 1998.} Shimura sent me in 1998).   Weil mentions this conversation in his 1979 comments on his 1967 article, included in his collected works.\footnote{Weil, {\it Oeuvres scientifiques}, Vol III, p. 450.}   The chronology is important for Lang, since the conversation definitely took place before Weil published his result.   \medskip

{\it When I first read Weil's answer about ``one and the other set being denumerableÓ, I characterized it as ``stupidÓ. I have since also called it inane. But actually, Weil's answer gives further evidence that he did not think of the conjecture himself. Indeed, as a result of his conversations with Serre and Weil, {\bf Shimura was directly responsible for changing the prevailing psychology} about elliptic curves over $\QQ$.}\footnote{Lang, p. 1304.}  [my emphasis; we will be returning to psychology]\medskip
 
(b) Shimura added two arguments that need to be considered separately. In the first place, there is the important theorem that he first published in 1971 in his textbook {\it Introduction to the Arithmetic Theory of Automorphic Functions}:  namely, that

\begin{theorem*}[Shimura]  {\bf To a modular form of weight 2 for $\Gamma_0(N)$ with rational Fourier coefficients one can associate an elliptic curve with rational coefficients.}  
\end{theorem*}

This is the evidence Weil cited in 1967 (as a {\it Mitteilung} from Shimura) and about which he asks ``whether things always behave that way."   

Twenty years ago, Shimura was attaching a great deal of importance to being credited with the proof of this theorem.  He was was unhappy with Wiles's discussion of this and a related point in one of his articles on FLT: \medskip

{\it ``Strangely Wiles says `Any such elliptic curve has the property... has an analytic continuation ...'  without saying that the result is due to me. This is another defect. Does he think that it is trivial?"}\footnote{Shimura, op. cit., May 11, 1998.  In his letter to Serge Lang, February 14, 1999, appended to this article, he refers to ``one more annoying item" in Wiles's article:  ``He attributes that statement to nobody.  I wonder if he considers it trivial.  You may call it a non-reference by Wiles."

The letter to Lang mainly complains that Wiles failed to appreciate the importance of the boldfaced Theorem, and only mentions the result about analytic continuation in the final paragraph.   This observation could well serve as the starting point for an analysis of the virtues inherent in a new theorem; but I don't have enough information to reconstruct Shimura's thoughts on the matter.  
}\medskip

(c)  Shimura's third argument is of a different nature again.  In his remarks on receiving the American Mathematical Society's Steele Prize for Lifetime Achievement, Shimura explained his motivation:\medskip

{\it ...at first I thought that  ...  curves  obtained  from  a  division quaternion algebra $B$ over $\QQ$ might not be modular... but I realized that no nonmodular $\QQ$-rational elliptic curves could be obtained for the following reason: Eichler had shown [that] the Euler products on $B$ are already included in those obtained from elliptic modular forms [2]. The Tate conjecture on this was explicitly stated much later, but the idea was known to many people, and so {\bf it was natural for me to think} that two elliptic curves with the same zeta function are isogenous. This fact concerning $B$, in addition to the results I had about the zeta functions of modular curves, may have been the strongest reason for my stating the conjecture that every $\QQ$-rational elliptic curve is modular.}\footnote{Shimura, {\it Notices of the AMS}, November 1996, p. 1345.} [my emphasis, to be explained later]\medskip

\subsection*{3.  The case for Weil}  Of course Weil gets credit for writing and publishing the paper that includes his converse theorem, an insight at the basis of an entire research program.  A vast generalization of the converse theorem by James Cogdell and Ilya Piatetski-Shapiro has been indispensable in 20 years of progress on the Langlands program.

Weil may have remembered the conjecture on elliptic curves from that 1964 conversation with Shimura.  This doesn't strengthen the case for attaching his name to the conjecture.  But Weil's paper included a major advance, already mentioned above, on the previous versions of the conjecture:  he suggested that an elliptic curve $E$ is a quotient of the Jacobian of the modular curve with level subgroup $\Gamma = \Gamma_0(N)$, where $N$ is an explicit integer, the {\it conductor} of $E$, defined in terms of the arithmetic of $E$.  This makes the conjecture considerably more precise:  there is an equivalence between two sets of objects that matches a specific numerical property; moreover for a given $N$, the set of such objects is finite.  Serre wrote that this was {\it a beautiful new idea ; it was not in Taniyama, nor in Shimura (as Shimura himself wrote to me after Weil's paper had appeared)}."  But he then added an explicit value judgment:
\medskip

{\it Its importance comes from the fact that it made the conjecture {\bf checkable} numerically (while Taniyama's statement was not). I remember vividly when Weil explained it to me, in the summer of 1966, in some Quartier Latin coffee house. Now things really began to make sense. Why no elliptic curve with conductor $1$ (i.e. good reduction everywhere)? Because the modular curve $X_0(1)$ of level $1$ has genus $0$, that's why! I went home and checked a few examples of curves with low conductor : I did not know any with conductor $< 11$ nor with conductor 16?  No surprise, since $X_0(N)$ has genus $0$ for such values of $N$, etc. Within a few hours, I was convinced that the conjecture was true.}\footnote{Jean-Pierre Serre, ``L'histoire de la `modularity conjecture'", {\it Gazette des math\'ematiciens}, {\bf 91}, Janvier 2002, 55--57.}   
\medskip

The virtue Serre identified in Weil's formulation -- that it lends itself to numerical verification -- has been taken up by defenders of the claim for Weil's priority; most of the defenders I know happen to be French.

\section*{``Conceptual evidence"}

We should also hear what Andrew Wiles himself has to say about the origin of the modularity conjecture.  Wiles may have been thinking of the virtue mentioned by Serre in the first page of his paper on Fermat's Last Theorem:\medskip

{\it A well-known conjecture which grew out of the work of Shimura and Taniyama in the 1950's and 1960's asserts that every elliptic curve over $\QQ$ is modular.  However, it only became widely known through its publication in a paper of Weil in 1967..., in which, moreover, Weil gave conceptual evidence for the conjecture. }      \medskip                                                                                                                                                                                                                                                                                                                            

At the time, though, I understood that the ``conceptual evidence" was the derivation of modularity from the characterization of modular forms given by Weil's converse theorem, together with the Hasse conjecture, in Weil's version (A3), and the Tate conjecture, mentioned below.

It's important to emphasize that all three protagonists -- Taniyama, Shimura, and Weil -- justified their speculation in roughly the same terms.  Taniyama's 12th problem includes the following sentence:\medskip

{\it If a conjecture of Hasse is true ... then the Fourier series obtained from $L_C(s)$ by the inverse Mellin transformation must be an automorphic form of dimension $-2$, of some special type (cf. Hecke).}\medskip

Shimura never published any version of his conjecture, but Lang reported in 1995 that\medskip

{\it The rationale for Shimura's conjecture was precisely the conjectured functional equation (Hasse), along the lines indicated in Taniyama's problem 12, suitably corrected. Shimura's bolder insight was that the ordinary modular functions ... suffice to uniformize elliptic curves defined over the rationals.}
\medskip

Shimura confirmed this a few years later in a private communication.\footnote{Shimura, unlike Wiles, did not see Weil's paper as a conceptual breakthrough.  In his e-mail to me dated May 11, 1998, he wrote 
```Conceptual evidence' is  rather misleading.  After all Taniyama and I were familiar with the results of Hecke on the characterization of an automorphic form by the functional equation of its Mellin transform, and  our ideas were partly based on that correspondence. É I always viewed Weil's result as its nice refinement.  But ``conceptual"?  Really? I don't get it."

The word ``conceptual" clearly continued to bother him, because in the letter to Serge Lang cited in footnote 10, he complained vigorously about Wiles's failure to understand his (Shimura's) contribution.  Shimura sent me a copy of this letter and I have included a reproduction of the letter at the end of this article.  The letter to Lang makes it clear, however, that Shimura was especially annoyed that Wiles and others failed to give him credit for the boldface Theorem cited above, under point (b) of the case for Shimura, and which he calls a ``geometric construction."}
  And Hasse's conjecture  constituted the ``certain theoretical grounds" [{\it gewisse theoretische Gr\"unde}] to which Weil alluded (elliptically) in his paper.    However, Hasse's conjecture alone would have indicated that the function $L_C(s)$ satisfied some functional equation, but not necessarily the condition (A2) that had been shown by Hecke to characterize modular forms.  So Weil actually took his speculations to the next step:  ``{\it Weiter l\"asst sich vermuten}",\footnote{Hard to translate.  On www.linguee.de one finds ``it looks as if," ``it can be assumed that," and ``suggests that" among the translations, depending on the context.}   as indicated above, that $L_C(s)$ satisfies Weil's  condition (A3), which Serre thought added the necessary precision to make the speculation of Taniyama and Shimura into a checkable and therefore full-fledged conjecture.

\section*{Discussion:  Naming Virtues}

\subsection*{1.  Taniyama's clear and distinct idea}

Taniyama's reasoning follows a pattern that is familiar to working mathematicians but that to my knowledge has not been taken up by philosophers.  ``If a conjecture of Hasse is true," he says, then something else must follow, because of Hecke's well-known characterization of modular forms.  Taniyama's 12th problem can be taken to be a sketch of a proof of a theorem:

\begin{theoremT*} If Hasse's conjecture is true, then the $L$-function of any elliptic curve over a number field $k$ is an automorphic form of weight $2$ of some special type.
\end{theoremT*}

In more detail:  (i) Hasse's conjecture implies that $L_C(s)$ satisfies a functional equation, and (ii) Hecke's theorem then implies that $L_C(s)$ is the L-function of an automorphic form.  Step (ii) is false for at least two reasons:  (a) Hecke's theorem only applies when k is the field of rational numbers, and (b) it only applies when the conductor $N$ is $1$, which is never the case (a point Serre mentioned in the text cited above).  So the virtue of Taniyama's contribution seems to involve the sketch of an incorrect proof of an imprecisely stated and undoubtedly false theorem (see (a)).  The mistakes do not enhance the contribution, of course; the imprecise and partially incorrect reasoning served as the basis of Taniyama's imaginative leap, making him the first to intuit an unexpected connection between elliptic curves and modular forms.  

The imaginative leap was not totally unprecedented.  Dirichlet series with Euler products and functional equations arise in number theory in two ways.  
One has to do with Galois theory, the other involves constructions in complex analysis that can be interpreted as some sort of harmonic analysis on topological groups.\footnote{This perspective was already visible in Erich Hecke's theory of L-functions of Gr\"ossencharaktere, and was made systematic in John Tate's influential thesis of 1950.  Kenkichi Iwasawa had independently discovered an essentially identical method slightly earlier.  Although Iwasawa never published his results and never wrote them up in detail in English, there has never been a priority dispute.} 
Class field theory, the major accomplishment of algebraic number theory in the first half of the 20th century, proves that both ways give rise to the same family of Dirichlet series when the structures involved (the Galois groups and the topological groups) are both abelian.  
Taniyama can be credited with the insight that something analogous is going on when the groups are not abelian.  To the best of my knowledge, this sort of insight is considered contingent and thus not suitable as a topic for philosophical analysis.  
It is undoubtedly the most prized of all the mathematical virtues, however, and for that reason alone deserves the attention of anyone who is interested in mathematical practice.  For example, the mathematical community would see no virtue in the publication of a new speculative intuition every week on the hope that one of them will turn out to be right.   A new conjecture not accompanied by a compelling motivation -- what we might call a {\it clear and distinct idea} -- will win no credit for its author.  In this respect the reasoning that acoompanied Taniyama's imaginative leap qualifies as {\it cartesian} in Ian Hacking's sense:  any reader familiar with both Hasse's conjecture and Hecke's theorem will immediately see the pertinence of Taniyama's 12th problem.   

The apparent importance of the cartesian character of Taniyama's insight has an apparently paradoxical implication that deserves to be stressed:  there is more virtue in Taniyama's incorrect formulation of his 12th problem, accompanied as it was by steps (i) and (ii) above, than there would have on the counterfactual assumption that he had presented with no justification whatsoever the fully-formed and correct modularity conjecture, as later verified by Breuil-Conrad-Diamond-Taylor.  
Something about this paradox feels wrong, but I am not sure what it is.

Of course, mathematicians generally show more appreciation for a conjecture when it is accompanied by a proof, however un-cartesian it might be, than when it is merely accompanied by a clear and distinct idea.  Nevertheless, although Hacking thought cartesian proofs are rare, I would argue that cartesian insights are fundamental to mathematical communication and progress.    Everyone remembers that Caesar crossed the Rubicon; only eyewitnesses and specialized historians remember the details of what happened when he got to the other side.  

Before turning to the other two claimants to priority, I remind the reader that the 13th problem on Taniyama's list was a version of Conjecture II.  I have alluded indirectly to the possession of the ``clear and distinct idea" underlying this intuition on the part of each of the three protagonists.  In Taniyama's case, it is implicit in the words ``Concerning the above problem" that form the transition between his problems 12 and 13, together with his use of the word ``isogeneity."  This is a basic concept in the theory of {\it abelian varieties}, a class of algebraic varieties to which elliptic curves belong.  An abelian variety is a projective algebraic variety which is also an algebraic group -- its group structure is given by algebraic maps.   To algebraic curve is canonically assigned an abelian variety called its {\it Jacobian}; the Jacobian of the modular curve $X(\Gamma)$ is denoted $J(\Gamma)$.  The theory of abelian varieties is one of the best understood parts of algebraic geometry.  In particular, every abelian variety is a product of ``simple factors" but only ``up to isogeny."  An isogeny between abelian varieties is an isomorphism up to finite error:  it is a surjective homomorphism of algebraic groups with finite kernel.  Taniyama's problem 13 says the Jacobian $J(\Gamma)$ for $\Gamma$ of ``stufe" $N$ ``contains" certain a elliptic curve; by ``contains" he means that it is one of its simple factors Ð up to isogeny.

We have seen that an L-function can be attached to any elliptic curve, but in fact Hasse's conjecture (or the Hasse-Weil conjecture) concerns the L-function attached to  any algebraic variety.  In particular, one can attach an L-function to an abelian variety.  It has been known for a long time that two isogenous abelian varieties have the same L-function.  Recall what Shimura wrote in 1996:  {\it it was natural for me to think that two elliptic curves with the same zeta function are isogenous}.  For want of a better account of the {\it nature} to which Shimura is referring here, let's assume that this ``idea," of which he claimed that it was ``known to many people," was ``clear and distinct" in their minds.   It was this idea to which Weil alluded when we wrote that it was ``of course obvious -- {\it nat\"urlich naheliegend zu erwarten}."  Here is the full citation:

``It is of course obvious to expect that under these conditions $C'$ [the elliptic curve factor of the Jacobian attached to the modular form by Shimura] is isogenous to $C$ [the original elliptic curve that on ``certain theoretical considerations" should be modular]; this is in fact confirmed in some cases." ({\it Es ist nat\"urlich naheliegend zu erwarten, da§ unter diesen Umst\"anden $C'$ mit $C$ isogen ist; das best\"atigt sich tats\"achlich in einigen F\"allen.})

With this in mind, problem 13 can be read in two ways.  First, is it true that every $J(\Gamma)$ of ``stufe" N has elliptic curves among their simple factors?  This is the most obvious reading of the text, but it is a strange question, because there are many modular curves whose Jacobians do not involve any elliptic curves.  The question was given a precise answer by Shimura's analysis of the Jacobians of modular curves in his 1971 textbook. 

The alternative reading is that every elliptic curve with rational coefficients is isogenous to a factor of the Jacobian of some modular curve.  This is now known to be true, as a consequence of two very substantial theorems: the theorem of Breuil-Conrad-Diamond-Taylor on the modularity of elliptic curves, and the 1983 theorem of Gerd Faltings that includes the {\it Tate conjecture} that Shimura had mentioned:  that two abelian varieties over a number field whose L-functions coincide are necessarily isogenous.\footnote{At some point in the 1980s I had the disconcerting experience of being asked by Shimura whether or not it was true that Faltings had proved this theorem.  Unfortunately, I don't remember whether it was before or after Faltings had joined Shimura in the Princeton mathematics department, nor whether it was before or after the work of Frey and Ribet, around 1985-86, that showed how Fermat's Last Theorem would follow from the Modularity Conjecture, in view of Faltings's theorem.  Either way, Shimura was surrounded by people who were far better qualified than I to answer his question.  

Tate's conjectures on algebraic cycles, including the one about isogenies of abelian varieties that Faltings proved, were made public in the 1960s.  They are fundamental to our thinking about the Galois representations attached to the l-adic cohomology of algebraic varieties, but I am not familiar with their prehistory and cannot comment on why Taniyama, Shimura, and Weil found it ``natural to think" that it was valid for elliptic curves.  
}    This is the expectation that Shimura found ``natural" and Weil found ``obvious," once Conjecture I has been established.  While one needs a good deal of expertise and sophistication in order to perceive Conjecture II as an obvious or natural consequence of Conjecture I and Tate's Conjecture, the choice of words suggests that neither Weil nor Shimura is claiming priority for this particular aspect of the conjecture.  So  I will henceforth restrict my attention to Conjecture I.

\subsection*{2.  Lang's realism, Shimura's phenomenology}

Lang's case for Shimura's priority is primarily forensic in nature:  his detective work places Weil and Serre in conversations with Shimura on this topic during a precise period in the early 1960s.  Weil and Serre were skeptical when Shimura asserted that he ``believed [that an elliptic] curve [over $\QQ$] should always be a quotient of the Jacobian of a modular curve" (both Conjectures I and II, in other words).  Since Weil was not only a skeptic at the time but even responded in a way Lang characterized as ``stupid" and ``inane," how could anyone assign him any credit for the conjecture?

If Lang's approach can be assigned to any philosophy, I would tend to call it {\it empiricist} and {\it realist}:  to the extent that there is a fact of the matter (and Lang clearly believes there is), it can be established by gathering evidence.  But Lang made a comment that calls this assignment into question:  as we have seen, he claimed that ``Shimura was directly responsible for changing the prevailing psychology about elliptic curves over $\QQ$."  We can agree that something changed about the way number theorists were thinking about elliptic curves over $\QQ$ without commiting ourselves to Lang's use of the word ``psychology," nor do we have to accept his judgment that Shimura had direct responsibility for the change by speaking with Weil and Serre.  Let's say that we could ascribe  to the community of number theorists the {\it disposition} to think about elliptic curves over $\QQ$ in certain ways in 1960, and a {\it markedly different disposition} in 1970.  Then we can say that Lang counts among the prized virtues of mathematical practice {\it the ability to alter dispositions} with regard to important topics, and that this is one of the main virtues he identifies in Shimura's role in proposing the modularity conjecture.  This doesn't suffice as a characterization of Shimura's contribution, however.  Other distinguished number theorists, including Barry Mazur,  Yuri Manin, and especially John Tate, as well as Serre himself, were busily altering dispositions with regard to elliptic curves over $\QQ$ in the decade following the publication of Weil's article, by giving lectures and writing survey articles\footnote{Especially J. Tate, The arithmetic of elliptic curves, {\it Invent. Math.}, {\bf 23} (1974), 179--206, which is where I and everyone else I knew developed our dispositions at the time. }  promoting a very precise and detailed international research program that continues to frame thinking about elliptic curves even now.  One presumes that they did so with dispositions altered as a direct result of Shimura's intervention with Serre and Weil a few years earlier -- and, I should add, with Tate (``I told it to Tate around the same time, and I remember that he was impressed."\footnote{Shimura, letter to the author, May 11, 1998.}). 

Lang is no longer among us, and we don't know what he had in mind when he refers to Shimura's ``changing the prevailing psychology."  My guess is that he would like his readers to believe that Weil's conversation with Shimura inspired him to think about the questions that led to his 1967 paper, inviting us to interpret his failure to mention either the contents of the conversation or Taniyama's 1955 questions as a deliberate attempt to claim the virtue of the original idea of Taniyama and Shimura for himself.  That interpersonal conflict does not lend itself to philosophical analysis, whereas ``dispositions to think about" various mathematical objects can be legitimately understood as targets of virtuous intervention within the community of mathematical practitioners.  The dispositions, and the interventions to change them, can thus be incorporated in a philosophy of mathematical practice, whether or not Shimura's private conversations were as effective as Lang claimed.  At the very least, this way of looking at things raises interesting questions about the relative virtue of private and public interventions within the community -- and we should not underestimate the importance of the letters of recommendation with which this essay began in creating dispositions in the first place.  

Shimura's own arguments (b) and (c) make no reference, explicit or otherwise, to the community of mathematical practitioners.  To put them in perspective, it's helpful to rewind to the early 1950s, when Hasse's conjecture was known in very few cases but was the basis of a very active research program.  Weil, Taniyama, and Shimura were reporting on their contributions to this program as early as the 1955 Tokyo-Nikko Conference, but their results (like those of Max Deuring, also present at the conference) were based on what was already known about abelian class field theory.  The first non-abelian example had been treated in a 1954 article by Martin Eichler.  Shimura spent much of the next 20 years generalizing Eichler's observation and creating what is now called (following Pierre Deligne) the theory of {\it Shimura varieties}, of which the simplest (but by no means simple) case is the one treated in his 1971 textbook.  In that book, and in a series of papers dating back to 1958, Shimura developed a systematic method for verifying new non-abelian cases of Hasse's conjecture for elliptic curves, in the process providing evidence for both Conjectures I and II.  Lang puts it this way: \medskip

{\it Shimura himself in the late fifties and sixties extended Eichler's results and proved that elliptic curves which are modular have zeta functions which have an analytic continuation. (Cf. the three papers [Sh 58], [Sh 61], and [Sh 67].)}\footnote{Lang, p. 1303.}\medskip

Thirty years earlier, this had already been mentioned in a survey article of J.W.S. Cassels that was largely responsible for creating the ``prevailing psychology" about elliptic curves in the 1960s:

\medskip
{\it Hasse's conjecture has also been verified for a few C  without complex multiplication which arise in the theory of modular functions [Shimura (1958a, 1961a, 1964a)]}\footnote{Cassels, J. W. S., Diophantine equations with special reference to elliptic curves. 
{\it J. London Math. Soc.}, {\bf 41} (1966) p. 280. }
\medskip

What should we call the virtue, exemplified by this evidence, on which Shimura bases his claim (b) to priority?  Wikipedia suggests one possible name:
\medskip

{\it Shimura's approach, later presented in his monograph, was largely {\bf phenomenological}, pursuing the widest generalizations of the reciprocity law formulation of complex multiplication theory.}[my emphasis]\footnote{Wikipedia, ``Shimura variety,"  https://en.wikipedia.org/wiki/Shimura\_variety}
\medskip
 
Wikipedia's history page informs me that on March 31, 2007, Charles Matthews inserted the word that I have set in boldface, together with its philosophical overtones.  The word is not used in Husserl's sense but rather as it is understood in philosophy of science.   Shimura developed the theory in order to establish the facts of the matter.\footnote{A Wikipedia editor named Arcfrk emphasized (in 2010) the contrast with Deligne's definition, which served ``to isolate the abstract features that played a role in Shimura's theory."}   In the course of this development, he made the observation reported above in argument (c):  that the phenomena (elliptic curves) observed in his study of the Shimura varieties attached to quaternion division algebras provided additional evidence for Conjectures I and II.   

Thus Shimura claimed the virtue that Matthews calls {\it phenomenological}.  As he himself acknowledged, it built upon Taniyama's original {\it cartesian} insight, but from the standpoint of philosophical analysis it is conveniently different.  Nevertheless, I don't want to leave Shimura without pointing out that mathematical practice sees an unbridgeable gulf between {\it providing (phenomenological) evidence} for a conjecture and {\it proving} the conjecture.  This is obvious, of course, and since no one is claiming that any of the triad (Taniyama, Shimura, Weil) came close to proving the conjecture, its relevance to identifying the virtues of a conjecture may not be clear.  This may best be understood {\it a posteriori}, in the light of the Langlands program for automorphic forms, specifically the program Langlands outlined in 1979 for the study of Shimura varieties.  During the 1960s, Alexander Grothendieck and his school had developed the theory of \'etale cohomology, which made it possible to attach Dirichlet series to automorphic forms in both of the ways mentioned previously.  \'Etale cohomology provided the connection with Galois theory, while the relation to {\it representation theory} -- what we previously called harmonic analysis on topological groups -- gave a variety of constructions of Dirichlet series.  The generalization of Conjecture I (no uniform generalization is possible for Conjecture II) asserts roughly that 

\subsubsection*{{\bf A}.  Galois-theoretic Dirichlet series are representation-theoretic Dirichlet series}

\hspace{.5 cm}

Specialists call this the ``Galois-to-automorphic" direction of Langlands's conjectures; current work on this direction largely derives from the methods Wiles developed in his proof of Fermat's Last Theorem.   The ``automorphic-to-Galois" direction, naturally, asserts

\subsubsection*{{\bf B}.  Representation-theoretic Dirichlet series are Galois-theoretic Dirichlet series}

\hspace{.5 cm}

Shimura's phenomenological efforts established {\bf B} for the class of modular forms that is relevant to elliptic curves with rational coefficients, and provided the starting point for all subsequent work on {\bf B}, which by now is established in considerable (though by no means complete) generality.  But it's not immediately obvious why even a complete proof of {\bf B} should help to establish a disposition to believe {\bf A}, as it apparently did for Shimura.  

Shimura actually did prove {\bf A} -- Conjecture I for the class of {\it elliptic curves with complex multiplication} -- in a paper published in 1973.\footnote{Goro Shimura, On the factors of the jacobian variety of a modular function field,  {\it J. Math. Soc. Japan}, {\bf 25} (1973), 523--544.  This theorem can be proved using Weil's converse theorem, but this was not the method Shimura chose, although his paper on modular forms of half-integral weight, published the same year, was based on a highly original application of Weil's 1967 paper.}   Even though this class of elliptic curves has especially favorable properties (the problem in this case is very close to being abelian), phenomenological results of this kind typically do count as significant evidence for {\bf A} within the framework of the Langlands program, since it is generally believed that we will have to wait several generations, if not centuries, for a complete proof of {\bf A}.  

\subsection*{3.  Weil's falsificationism, on Serre's reading}

Unlike some of his Bourbaki colleagues, Serre has never displayed an interest in philosophy for its own sake, but the virtue he discerns in Weil's contribution to the conjecture is the easiest to translate into familiar philosophical terms.   If we replace the word ``checkable" with its synonym ``verifiable," then Weil's virtue, according to Serre, is that he has endowed the Modularity Conjecture with a {\it meaning}, following the {\it logical positivist} dictum that (in the words of Schlick) ``The meaning of a proposition is the method of its verification."\footnote{Moritz Schlick, ``Meaning and verification," {\it The Philosophical Review}, {\it 45}, No. 4 (Jul., 1936), 339-369.}    In practice, though, {\it falsification} of the Modularity Conjecture by finding a counter-example would have carried more weight than its verification in special cases, if for no other reason than that a proof by (philosophical rather than mathematical) induction would require infinitely many verifications, which is impossible, whereas a single counter-example suffices for disproof.   Thus the virtue Serre highlights may just as well be {\it popperian}.  Either way, by insisting on the virtue of making a conjecture ``checkable," Serre seems to be in agreement with Mach, when he asserted that ``where neither confirmation nor refutation is possible, science is not concerned."\footnote{Mach, Ernst, 1883, {\it Die Mechanik in ihrer Entwicklung}, Leipzig: Brockhaus, 9th ed. transl. by T.J. McCormick, {\it The Science of Mechanics}, Chicago: Open Court, 1960, p. 587.}

The clarity of Serre's argument simplifies the philosopher's task; he has identified a virtue in conjecture-making that is both unambiguous and easy to map onto familiar concerns in the philosophy of science.  I would go further and say that Serre has rendered a service to philosophy:  while both the verificationism of the logical positivists and Popper's falsificationism are known to be vulnerable to critique, on the grounds that there is more to the world than its logical Aufbau, Serre's ``checkability" may well be a viable criterion for the meaningfulness of a conjecture.  At any rate, it has not been subject to philosophical analysis.  

Most importantly, perhaps, Serre has shown some consistency in following his own dictum:  Serre formulated his famous Conjecture on modular representations of Galois groups in 1973 but only published it in 1987, after, as he put it in his thanks,\medskip

{\it Jean-Franois Mestre ...a r\'eussi \`a programmer et v\'erifier un nombre d'exemples suffisant pour me convaincre que la conjecture m\'eritait d'\^etre prise au s\'erieux.}\footnote{Jean-Pierre Serre, Sur les repr\'esentations modulaires de degr\'e $2$ de $Gal(\overline{\QQ}/\QQ )$, {\it Duke Mathematical Journal}, {\bf 54} (1): (1987) p. 180.   Its proof some 20 years later, by Chandrashekhar Khare and Jean-Pierre Wintenberger, made extensive use of the methods introduced by Wiles, Taylor, and others in proving the Modularity Conjecture, to which Serre's conjecture is very closely related.} \medskip

The examples occupy 11 of the 52 pages of his manuscript.

\section*{Conclusion:  the ``prevailing psychology"}

Lang indirectly hinted at a model of the {\it mathematical subject} by assigning to Shimura direct responsibility ``for changing the prevailing psychology" by privately conveying his belief in the truth of what we are calling the Modularity Conjecture to two or three individuals.  This subject is structured hierarchically, so that impressions acquired by those (like Weil, Serre, or Tate) at the pinnacle of the hierarchy have the power of detectably altering dispositions throughout the subject.  Lang doesn't explain the means of this alteration, so we don't know the relative weight in the process of private exchanges, seminar presentations, or publications -- just to mention the mechanisms with which most members of the community would be familiar.  

But even granting that the dispositions of a collective subject can be altered as a result of a few brief conversations involving a group of exceptionally influential individuals, the question remains:  why should there have been a ``psychology about elliptic curves over $\QQ$" in the first place?    Whether or not there is a virtue ethics for mathematics, in which every conceivable topic merits its own psychology, in practice the attention of the mathematical subject is rather narrowly focused.  Why were elliptic curves so much on the minds of number theorists in the 1960s?

It's easier to understand this attention when we remember that number theory at its core is concerned with understanding {\it Diophantine equations} -- polynomial equations in several variables with integer coefficients-- and more specifically with characterizing their set of solutions in integers, or in rational numbers.   Equation \eqref{ell} that we encountered at the beginning of the article, the one whose solution set is an elliptic curve, is a Diophantine equation -- it is the simplest possible Diophantine equation whose theory has not been fully understood since the 19th century (if not longer) -- and yet the nature of its solution set remains an outstanding open question.  Much was known by the 1960s, however.  It was proved by L. E. J. Mordell nearly 100 years ago that its set of rational solutions forms a finitely generated abelian group, a result generalized in Andr\'e Weil's thesis and now known as the {\it Mordell-Weil Theorem}.  This landmark in number theory, arguably the deepest result known at the time about Diophantine equations, posed a challenge to number theorists:  how many generators does the group have?  And how can we tell whether the group is infinite or finite?  The challenge was more acute in view of the special status of elliptic curves among all algebraic curves\footnote{For our present purposes, an algebraic curve is defined by an equation in two variables, like \eqref{ell}.}  with rational coefficients.  Every such curve has a genus, a non-negative integer; elliptic curves are the curves of genus $1$.  The rational solutions of curves of genus 0 -- like the circle with equation $x^2 + y^2 = 1$ -- are well understood (rational points on the circle are essentially the same as right triangles with integer sides).  The {\it Mordell Conjecture} -- proved by Faltings in 1983 -- is the statement that curves of genus $2$ have only finitely many solutions.  For elliptic curves, simple as they are, there is no elementary way to determine whether the (rational) solution set is finite or infinite.

But in the early 1960s, when Shimura was privately conversing in Princeton with Weil, Serre, and Tate, over in England Bryan Birch and H. P. F. Swinnerton-Dyer were busy formulating and testing a {\it conjectural answer} to the question.  The {\it Birch-Swinnerton-Dyer (BSD) Conjecture} -- the appropriateness of its name has never been called into question! -- asserts that the number of generators of the group of rational points of the elliptic curve $C$ can be determined by examining its associated Dirichlet series -- Taniyama's $L_C(s)$, Weil's $\Lambda$ -- but only if the Dirichlet series had (at least some of) the properties conjectured by Hasse.  

The (eminently checkable) BSD Conjecture, which was widely known\footnote{Thanks to the 1966 article by Cassels cited in note 20, and of course to lectures by Birch and Swinnerton-Dyer themselves, especially the latter's chapter in the book {\it Algebraic Number Theory} edited by Cassels and A. Frhlich Ñ London:  Academic Press (1967) Ñ which for 50 years has been required reading for every student of number theory.}  by the time Weil's article was published, gave such an appealing answer to the questions raised by the Mordell-Weil Theorem that one can easily understand why number theorists wanted the Conjecture to be meaningful -- and thus for Hasse's conjectures to be true for elliptic curves, and thus for Taniyama's ``theorem" to be true.  With the help of the BSD Conjecture, the Modularity Conjecture thus becomes an object of number theorists' collective desire -- a symptom of our mass psychology, as Lang intuited.\footnote{The importance of the BSD Conjecture in promoting the Modularity Conjecture is beyond the scope of this article, but two landmark articles written in the decades following the publication of Weil's paper give a sense of the excitement it generated.   The introduction to ``Arithmetic of Weil Curves", by Mazur and Swinnerton-Dyer, published in 1974 (but available several years earlier) explicitly identifies the ``aim of Part I of this paper":   to ``examine... structures" that are ``associated to the Weil parametrization and the bearing they have (via the conjectures of Birch and Swinnerton-Dyer) on the arithmetic" of elliptic curves over $\QQ$.  Part II as well as Part I became instantly applicable to all elliptic curves over $\QQ$ with the proof of the Modularity Conjecture.   In 1983 Benedict Gross and Don Zagier announced their proof of their famous formula the first that explicitly constructs infinitely many rational points on an elliptic curve using information about its Dirichlet series.  Since the Modularity Conjecture was not  yet established, their article (published in 1986, before the controversy over priority broke out) explains that the results apply to ``strong Weil curves" and refer to the Mazur-Swinnerton-Dyer article for the terminology.  B. Mazur and H. P. F. Swinnerton-Dyer, Arithmetic of Weil Curves, {\it Inventiones Math.}, {\bf 25} (1974) 1--61; B. Gross and D. Zagier, Heegner points and derivatives of L-series,  {\it Inventiones Math.} {\bf 84} (1986) 225--320.}

\includegraphics[scale = .65]{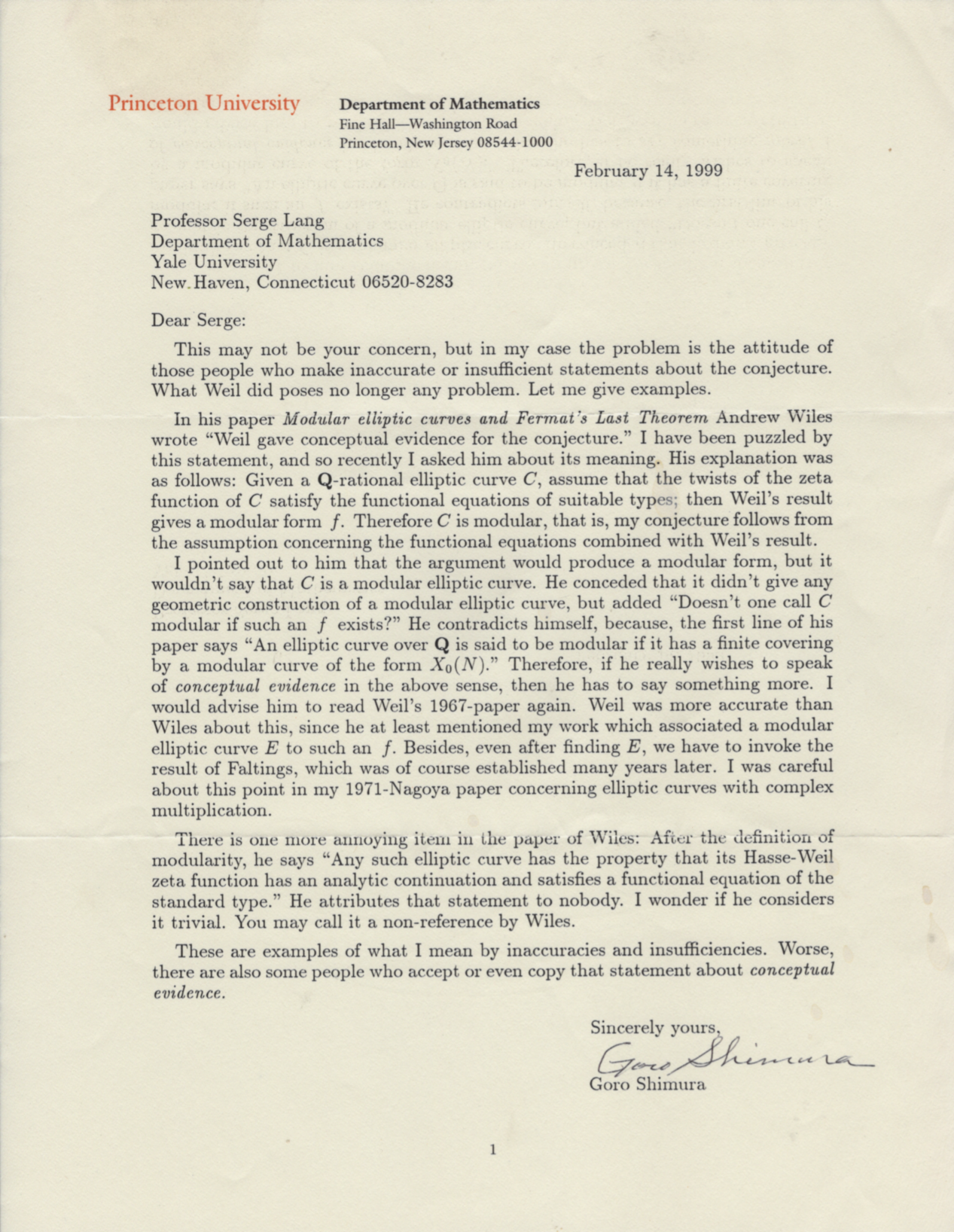}

\end{document}